\documentclass[12pt]{article}

\usepackage{amssymb}
\usepackage{amsfonts}
\usepackage{times}
\usepackage{cancel}
\usepackage{mathptmx}
\usepackage{amsmath}
\usepackage[usenames]{color}
\usepackage{amscd}
\usepackage{cite}
\usepackage{cases}
\usepackage{amsthm}
\usepackage{hyperref}
\usepackage{enumerate}
\usepackage{enumitem} 
\usepackage{indentfirst}
\usepackage{array} \usepackage{tabularray} \usepackage{caption} 
\usepackage{eucal} 

\hypersetup{colorlinks=true, linkcolor=black, citecolor=black}

\def\deg{{\rm deg}}

\def\cl{\centerline}
\def\vs{\vspace*}

\def\vs{\vspace*}

\allowdisplaybreaks[4]

\numberwithin{equation}{section}

\newtheorem{theorem}{Theorem}[section]
\newtheorem{definition}[theorem]{Definition}
\newtheorem{corollary}[theorem]{Corollary}
\newtheorem{lemma}[theorem]{Lemma}

\newtheorem{example}[theorem]{Example}

\newtheorem{remark}[theorem]{Remark}

\begin{document}

\begin{center}
    {\textbf{\large {Computation of cohomology of finite Lie conformal algebras}}}
    \footnote{$^\ast$Corresponding author: lmyuan@hit.edu.cn}
\end{center}
    
\cl{Zhenyu Zhang, Lamei Yuan$^\ast$}
    
\cl{\small School of Mathematics, Harbin Institute of Technology, Harbin
    150001, China\\}
\cl{\small E-mail: zhenyuzhang9@163.com, lmyuan@hit.edu.cn}
    
\vs{8pt}
    
\small\footnotesize
\parskip .005 truein
\baselineskip 3pt \lineskip 3pt
\noindent{\bf Abstract:}\vs{3pt} By using the energy operator introduced by B. Bakalov, A. De Sole, and V.\,G. Kac, we propose an algorithm for computing the cohomology of finitely generated free Lie conformal algebras with a Virasoro element.  
This computational method enables us to determine cohomologies with coefficients in their conformal modules, including the trivial module, irreducible finitely free modules, and adjoint modules. 
As concrete applications, we compute by \textit{Mathematica} both basic and reduced cohomologies for the $\mathcal{W}(b)$, Schr\"odinger-Virasoro and extended Schr\"odinger-Virasoro Lie conformal algebras.  
\vs{5pt}

		\vs{5pt}

    
\parskip .001 truein\baselineskip 6pt \lineskip 6pt

\section{Introduction}
The notion of a Lie conformal algebra (LCA) was introduced by Kac in \cite{KAC} to axiomatically describe the singular part in the operator product expansions of chiral fields in conformal field theory. Over the past three decades, Lie conformal algebras have been shown to exhibit deep connections with infinite-dimensional Lie (super)algebras, as well as assciative or non-assciative algebras. The general theory of structure and representation theory for Lie conformal algebras has been established in \cite{CJK} and \cite{DK1}, and widely developed in the literatures (see, for example, \cite{CKW,SXY1,SXY2,Yuan2}).

The cohomology theory of Lie conformal algebras, introduced in \cite{BKV} and further developed in \cite{DK2}, has attracted considerable attention due to its fundamental role in understanding the structure and deformation theory of these algebras\cite{Yuan2}. A central focus in recent years has been computing cohomology groups for Lie conformal algebras with coefficients in their modules.
The cohomology of the Virasoro conformal algebra $\rm Vir$ and the current algebra ${\rm Cur}\mathfrak{g}$, with $\mathfrak{g}$ simple, was determined in \cite{BKV}. The low-dimensional cohomology of the general conformal algebra $gc_N$ with coefficients in its trivial and natural modules was computed in \cite{Su}. With more concrete examples of Lie conformal algebras related to the Virasoro conformal algebra being constructed, their cohomologies have also been systematically studied. For instance, the cohomology of the Heisenberg-Virasoro conformal algebra with coefficients in its trivial, natural and irreducible finite conformal modules was computed in \cite{YW1}. 
Similar analyses were later extended to the Schr\"odinger-Virasoro and extended Schr\"odinger-Virasoro Lie conformal algebras \cite{WL1,WL2}.
However, the cohomology of finite Lie conformal algebras with coefficients in their adjoint modules remains unexplored. Moreover, there is no unified method for computing cohomology groups of Lie conformal algebras. These open questions serve as the primary motivation for the present work.


In \cite{Su,YW1,WL1,WL2}, the concept of ``homogeneous cochain" seemed to be an important tool in the computation of cohomology of Lie conformal algebras. Later in \cite{BDK}, B. Bakalov, A. De Sole, and V.\,G. Kac introduced the energy operator theory to compute the cohomology of Poisson vertex algebras containing a Virasoro element. Applying this to Lie conformal algebra cohomology, we find that the ``homogeneous cochain" is actually an eigenvector of the energy operator on the basic complex induced by a Virasoro element $L$, 
such that 
$$[L_\lambda L]=(\partial+2\lambda)L,$$
and such that
\begin{equation*} 
L_{(0)}:=[L_\lambda \cdot~]\big|_{\lambda=0}=\partial, \quad L_{(1)}:=\frac{\rm d}{{\rm d} \lambda}[L_\lambda \cdot~]\big|_{\lambda=0}\in {\rm End}~\mathcal{A}~ \text{is diagonalizable}.   
\end{equation*}
In this setting, $M$ is a conformal $\mathcal{A}$-module with respect to the Virasoro element $L$ in $\mathcal{A}$, such that
\begin{equation*}
L_{(0)}^M:=(L_\lambda \cdot~)\big|_{\lambda=0}=\partial^M, \quad L_{(1)}^M:=\frac{\rm d}{{\rm d} \lambda}(L_\lambda \cdot~)\big|_{\lambda=0}\in {\rm End}~M~ \text{is diagonalizable}.
\end{equation*}
The energy operator on $M[\lambda_1,\cdots,\lambda_n]$ is defined by
\begin{align*} E=L_{(1)}^M+\sum_{i=1}^{n}\lambda_i \frac{\rm d}{{\rm d}\lambda_i}, \end{align*}
which further induces a diagonalizable energy operator $\widetilde{E}$ on the basic complex $\widetilde{C}^q(\mathcal{A},M)$. Furthermore, we can decompose every $q$-cochain $f\in\widetilde{C}^q(\mathcal{A},M)$ into a formal sum of $\widetilde{E}$-eigenvectors, so that  
\begin{align*}\label{decompsiton1}
f=\sum_{r\in\mathbb{F}} f^{(r)},
\end{align*}
where the ``homogeneous cochain" $f^{(r)}\in \widetilde{C}^q(\mathcal{A},M)$ is an eigenvector of $\widetilde{E}$ with eigenvalue $r$. Based on this, we prove the main theorem of the paper (Theorem \ref{3.2.4}) stating that the eigenvalues of the energy operator $\widetilde{E}$ on $\widetilde{\mathrm{H}}^q(\mathcal{A},M)$ can be only $0$. This can be seen as a generalizition of Lemma 3.2 in \cite{Su}. We extend their method to a more general setting: free LCAs containing a Virasoro element, which range from finite LCAs including the Virasoro LCA, $\mathcal{W}(b)$ LCA and Schr\"odinger-Virasoro LCA, to infinite LCAs including the loop Virasoro LCA, general LCA and Block type LCA $\mathfrak{B}(p)$ with $p\neq0$. Then we employ the energy operator technique combined with analysis of a long exact sequence to derive the fundamental dimensional relationship between basic and reduced cohomology: 
$$\dim~ \mathrm{H}^q(\mathcal{A},M)
=\dim~ \widetilde{\mathrm{H}}^q(\mathcal{A},M) + \dim~ \widetilde{\mathrm{H}}^{q+1}(\mathcal{A},M),\quad \forall~q\geq1.$$

Based on those results, in the last section, we develop an explicit algorithm for computation of both $\widetilde{\mathrm{H}}^q(\mathcal{A},M)$ and ${\mathrm{H}}^q(\mathcal{A},M)$, provided that $\mathcal{A}$ is an LCA with a Virasoro element, $M$ is a conformal $\mathcal{A}$-module, and both of them are finitely and freely generated as $\mathbb{F}[\partial]$-modules. Then by \textit{Mathematica}, we compute basic and reduced cohomologies of the $\mathcal{W}(b)$, Schr\"odinger-Virasoro and extended Schr\"odinger-Virasoro conformal algebras with coefficients in their adjoint modules. The corresponding results are presented in Theorems \ref{w0}--\ref{sv}. Furthermore, we derive some partial cohomological results for them with coefficients in their finitely free irreducible modules.

Throughout the paper, the base field $\mathbb{F}$ is a field of characteristic $0$. Unless
otherwise specified, all vector spaces and their tensor products are over $\mathbb{F}$. We denote by $\mathbb{Z}_+$ and $\mathbb{Z}_{-}$ the sets of all non-negative integers and negative integers, respectively. 

\section{Preliminaries}

In this section, we will review the definitions of a Lie conformal algebra (LCA), a module
over it and the corresponding LCA cohomology. 

\begin{definition}{\rm (\kern-3pt\cite{KAC})
A Lie conformal algebra (LCA) $\mathcal{A}$ is an $\mathbb{F}[\partial]$-module equipped with an $\mathbb{F}$-bilinear map $[\cdot_\lambda\cdot]:\mathcal{A}\otimes \mathcal{A}\longrightarrow \mathcal{A}[\lambda]=\mathbb{F}[\lambda]\otimes\mathcal{A}$, $a\otimes b \mapsto [a_\lambda b]$, called the $\lambda$-bracket, which satisfies the following axioms ($a,b,c\in \mathcal{A}$):
\begin{enumerate}[label=\arabic*]
    \item[(L1.)] Conformal sesquilinearity:$~~[\partial a_\lambda b]=-\lambda[a_\lambda b]$, $[a_\lambda \partial b]=(\partial+\lambda)[a_\lambda b]$;
    \item[(L2.)] Skew-symmetry:$~~[a_\lambda b]=-[b_{-\lambda-\partial} a]$;
    \item[(L3.)] Jacobi identity:$~~[a_\lambda [b_\mu c]]=[[a_\lambda b]_{\lambda+\mu} c]+[b_\mu [a_\lambda c]]$.
\end{enumerate}}
\end{definition}

A Lie conformal algebra $\mathcal {A}$ is called {\it finite} if it is a finitely generated $\mathbb{F}[\partial]$-module; otherwise, it is called {\it infinite}. The rank of $\mathcal {A}$ is its rank as an $\mathbb{F}[\partial]$-module. 
For every $j\in\mathbb{Z}_+$, we can define the {\it $j$-product} $a_{(j)}b$ of any elements $a,~b$ in $\mathcal{A}$ by the following generating series:
\[ [a_\lambda b]=\sum_{j\in\mathbb{Z}_+}\frac{\lambda^j}{j!}(a_{(j)}b). \]
Then the following axioms of $j$-products hold for all $a,b\in \mathcal{A}$:
\begin{align*}
&a_{(n)}b=0,\ {\rm for}\ n\gg0;\\
(\partial a)_{(n)}b&=-na_{(n-1)}b,~ a_{(n)}(\partial b)=\partial(a_{(n)}b)+na_{(n-1)}b;\\
a_{(n)}b&=\sum_{j\in\mathbb{Z}_{+}}(-1)^{(n+j+1)}\frac{\partial^j}{j!} b_{(n+j)}a;\\
[a_{(m)},&b_{(n)}]=\sum_{j=0}^m \begin{pmatrix}
m\\j
\end{pmatrix}(a_{(j)}b)_{(m+n-j)}.
\end{align*}
Actually, one can also define Lie conformal algebras by using the language of $j$-products (c.f. \cite{KAC}).

It is well known that the
simplest nontrivial Lie conformal algebra is the Virasoro conformal algebra $\mathrm{Vir}$, which is a rank one
free $\mathbb{F}[\partial]$-module generated by a symbol $L$, such that 
 \begin{equation} \label{2.1}
 [L_\lambda L]=(\partial+2\lambda)L. 
 \end{equation}

\begin{definition}{\rm (\kern-3pt\cite{BKV})
A module $M$ over a Lie conformal algebra $\mathcal{A}$ is an $\mathbb{F}[\partial]$-module with a $\lambda$-action $a_\lambda m$, which is an $\mathbb{F}$-bilinear map $\mathcal{A}\otimes M \longrightarrow M[\lambda]=\mathbb{F}[\lambda]\otimes M$, $a\otimes m \mapsto a_\lambda m$, such that ($a,b\in \mathcal{A}$ and $m\in M$)
\begin{align*}
&(\partial a)_\lambda m=-\lambda a_\lambda m,\quad a_\lambda(\partial m)=(\partial+\lambda)a_\lambda m,\\
&a_\lambda(b_\mu m)-b_\mu(a_\lambda m)=[a_\lambda b]_{\lambda+\mu}m.
\end{align*}}
\end{definition}

If an $\mathcal{A}$-module $M$ is finitely generated as an $\mathbb{F}[\partial]$-module, then $M$ is called {\it finite}; otherwise, it is called {\it infinite}. And if $M$ has no nontrivial submodules, then $M$ is called {\it irreducible}. By definition, $\mathcal{A}$ itself can be seen as a module of $\mathcal{A}$, called the {\it adjoint module}, where the $\lambda$-action is the $\lambda$-bracket of $\mathcal{A}.$ Given an $\mathcal{A}$-module $M$, for each $j\in\mathbb{Z}_+$, we can also define the $j$-actions of $\mathcal{A}$ on $M$ by using the
following generating series: \[ a_\lambda m=\sum_{j\in\mathbb{Z}_+}\frac{\lambda^j}{j!}(a_{(j)}m), ~~\mbox{for}~~a\in \mathcal{A},~m\in M.\]


In the following let us recall the Hochschild cohomology complex $C^\bullet(\mathcal{A},M)$ for a Lie
conformal algebra $\mathcal{A}$ with coefficients in its module $M$ (see \cite{BKV} for details).
\begin{definition}{\rm For $q\in \mathbb{Z}_+$, a $q$-cochain  of a Lie conformal algebra $\mathcal{A}$
with coefficients in an $\mathcal{A}$-module $M$ is an  $\mathbb{F}$-linear map  
$$f:\mathcal{A}^{\otimes q} \longrightarrow M[\lambda_1,...,\lambda_q], ~~  a_1\otimes\cdots\otimes a_q \mapsto f_{\lambda_1,...,\lambda_q}(a_1,....,a_q),$$
satisfying the following conditions:
\begin{itemize}
 \item[(B1.)] Conformal antilinearity:$~~f_{\lambda_1,...,\lambda_n}(a_1,...,\partial a_i,...,a_n)=-\lambda_i f_{\lambda_1,...,\lambda_n}(a_1,...,a_n) $,
    \item[(B2.)] Skew-symmetry:$~~f_{\sigma(\lambda_1),...,\sigma(\lambda_n)}\big(\sigma(a_1),...,\sigma(a_n)\big)=\mathrm{sign}(\sigma)f_{\lambda_1,...,\lambda_n}(a_1,...,a_n)$, where $\sigma\in S_n$ is a permutation of $n$ elements.
\end{itemize}}

\end{definition}

Let $\mathcal{A}^{\otimes 0}=M$, so that a $0$-cochain is an element in $M$. The vector space of all $q$-cochains is denoted by $\widetilde{C}^q(\mathcal{A},M)$. The differential operator $\mathbf d: \widetilde{C}^q(\mathcal{A},M)\longrightarrow \widetilde{C}^{q+1}(\mathcal{A},M)$ is defined as follows,
\begin{equation} \label{2.2}
    \begin{split}
(\mathbf d f)_{\lambda_1,...,\lambda_{q+1}}(a_1,...,a_{q+1})=&\sum_{k=1}^{q+1}(-1)^{k+1} a_{k\lambda_{k}}f_{\lambda_1,...,\hat{\lambda}_k,...,\lambda_{q+1}}(a_1,...,\hat{a}_k,...,a_{q+1}) \\
&+\sum\limits_{{k<l;k,l=1}}^{q+1}(-1)^{k+l}f_{\lambda_k+\lambda_l,\lambda_1,...,\hat{\lambda}_k,...,\hat{\lambda}_l,...,\lambda_{q+1}}([a_{k\lambda}a_l],a_1,...,\hat{a}_k,...,\hat{a}_l,...,a_{q+1}),
    \end{split}
\end{equation}
where $a_1,...,a_{q+1}\in \mathcal{A}$, and the symbol $~\hat{}~$  means the element below it is missing. In particular, 
\begin{align}\label{zero-cochain}
 (\mathbf d m)_{\lambda}(a) & =a_\lambda m, ~~\mbox{if}~m\in M~\mbox{is~a~zero-cochain.}
\end{align} 
From \cite{BKV}, we know that $\mathbf d$ preserves the space of cochains and $\mathbf d ^2=0$, so that all the cochains form a complex, which is denoted by 
\[ \widetilde{C}^\bullet(\mathcal{A},M)=\bigoplus_{q\in \mathbb{Z}_+}\widetilde{C}^q(\mathcal{A},M), \]
and called the {\it basic complex}. 

At the same time, we can define an $\mathbb{F}[\partial]$-structure on $\widetilde{C}^\bullet(\mathcal{A},M)$ via
\[ (\partial f)_{\lambda_1,...,\lambda_q}(a_1,...,a_q)=(\partial^M+\sum_{k=1}^{q}\lambda_k)f_{\lambda_1,...,\lambda_q}(a_1,...,a_q), \]
where $\partial^M$ denotes the action of $\partial$ on $M$. We have $\mathbf d \partial=\partial \mathbf d$, and thus $\partial\widetilde{C}^\bullet(\mathcal{A},M)$ is a subcomplex of $\widetilde{C}^\bullet(\mathcal{A},M)$. The quotient complex 
\[ C^\bullet(\mathcal{A},M)=\widetilde{C}^\bullet(\mathcal{A},M)/\partial \widetilde{C}^\bullet(\mathcal{A},M)=\bigoplus_{q\in\mathbb{Z}_+}C^q(\mathcal{A},M)\] 
is called the {\it reduced complex}.

Note that $C^0(\mathcal{A},M)=M/\partial M.$ We denote by $\int: M\rightarrow M/\partial M$ the canonical quotient map.

\begin{definition}{\rm 
The basic cohomology $\widetilde{\rm H}^\bullet(\mathcal{A},M)$ of a Lie conformal algebra $\mathcal{A}$ with coefficients in an $\mathcal{A}$-module $M$ is the cohomology of the basic complex $\widetilde{C}^\bullet(\mathcal{A},M)$, and the (reduced) cohomology ${\rm H}^\bullet(\mathcal{A},M)$ is the cohomology of the reduced complex $C^\bullet(\mathcal{A},M)$.}
\end{definition}
 
For $f\in\widetilde{C}^q(\mathcal{A},M)$, it's called a {\it $q$-cocycle} if ${\bf d}f=0$ and a {\it $q$-coboundary} if there exists $g\in\widetilde{C}^{q-1}(\mathcal{A},M)$, such that $f={\bf d}g$. The space of all $q$-cocycles is denoted by $\widetilde{D}^q(\mathcal{A},M)$, and the space of all $q$-coboundaries is denoted by $\widetilde{B}^q(\mathcal{A},M)$, so that the basic cohomology $\widetilde{\rm H}^q(\mathcal{A},M)$ is the quotient space $\widetilde{D}^q(\mathcal{A},M)/\widetilde{B}^q(\mathcal{A},M)$.

Let $\mathcal{A}$ be an LCA and $M$ be an $\mathcal{A}$-module. A {\it Casimir element} is an element $\int m\in M/\partial M$ such that $\mathcal{A}_{-\partial}m=0$. The space of all Casimir elements is denoted by ${\rm Cas}(\mathcal{A},M)$. In particular, ${\rm Cas}(\mathcal{A},\mathcal{A})=\{\int a\in \mathcal{A}/\partial\mathcal{A}:[a_\lambda \mathcal{A}]\big|_{\lambda=0}=0\}$.

A {\it derivation} from $\mathcal{A}$ to $M$ is an $\mathbb{F}[\partial]$-module homomorphism $D:\mathcal{A}\rightarrow M$ such that \[ D([a_\lambda b])=a_\lambda D(b)-b_{-\lambda-\partial}  D(a),~~ \forall~a,b\in\mathcal{A}.\]  
A {\it derivation} is called {\it inner} if it has the form $D_{\int m}(a)=-a_{-\partial}m$, for some $\int m\in M/\partial M$. Denote by ${\rm Der}(\mathcal{A},M)$ the space of derivations from $\mathcal{A}$ to $M$, and by ${\rm Inder}(\mathcal{A},M)$ the subspace of inner derivations. In the special case of $M=\mathcal{A}$, we have the usual
definition of a derivation and an inner derivation of the LCA $\mathcal{A}$ as an $\mathbb{F}[\partial]$-module endomorphism $D$ such that 
\[ D([a_\lambda b])=[a_\lambda D(b)]+[D(a)_{\lambda}b] ,\]
and, respectively, 
\[({\rm ad\,}a)b=a_{(0)}b=[a_\lambda b]\big|_{\lambda=0}.\]

\begin{theorem}{\rm (\kern-2.5pt\cite{BKV,DK2}) Let $\mathcal{A}$ be a Lie conformal algebra, and $M$ be an $\mathcal{A}$-module. Then
  \begin{itemize}
    \item[(1)] ${\rm H}^0(\mathcal{A}, M) = {\rm Cas}(\mathcal{A}, M).$
    \item[(2)] ${\rm H}^1 (\mathcal{A}, M) = {\rm Der}(\mathcal{A}, M)/ {\rm Inder}(\mathcal{A}, M).$
    \item[(3)] ${\rm H}^2 (\mathcal{A}, M)$ is the space of isomorphism classes of $\mathbb{F}[\partial]$-split extensions of the LCA $\mathcal{A}$ by the $\mathcal{A}$-module $M$, where $M$ is viewed as an LCA with zero $\lambda$-bracket. In particular, ${\rm H}^2 (\mathcal{A}, \mathcal{A})$ parameterizes the equivalence classes of first-order deformations of $\mathcal{A}$ that preserve the product and the $\mathbb{F}[\partial]$-module structure.
  \end{itemize}}
\end{theorem}

\section{The energy operator on LCA cohomology} 

In this section, we first recall the definitions of a Virasoro element $L$ over an LCA and conformal weight with respect to $L$. Then we introduce the notion of a diagonalizable energy operator $\widetilde{E}$ on basic LCA complex spaces with coefficients in conformal modules. This process is similar to \cite{BDK}. We show that the conformal weight in basic cohomology can be only zero, which is used in Sect.\ref{sec4} to design a unified algorithm to compute the basic and reduced LCA cohomology. 




\begin{definition}{\rm (\kern-2.5pt\cite{BDK})
A {\it Virasoro element} $L$ in a Lie conformal algebra $\mathcal{A}$ is an element satisfying \eqref{2.1} and
\begin{equation} \label{2.3}
L_{(0)}:=[L_\lambda \cdot~]\big|_{\lambda=0}=\partial, \quad L_{(1)}:=\frac{\rm d}{{\rm d} \lambda}[L_\lambda \cdot~]\big|_{\lambda=0}\in {\rm End}~\mathcal{A}~ \text{is diagonalizable}.   
\end{equation}
An LCA $\mathcal{A}$ is called {\it conformal} if it contains a Virasoro element.  We say that $a\in \mathcal{A}$ has {\it conformal weight} $\Delta(a)\in \mathbb{F}$ if it is an eigenvector of $L_{(1)}$ of eigenvalue $\Delta(a)$. 
An $\mathcal{A}$-module $M$ is called {\it conformal} with respect to the Virasoro element $L$ in $\mathcal{A}$ if
\begin{equation} \label{2.4}
L_{(0)}^M:=(L_\lambda \cdot~)\big|_{\lambda=0}=\partial^M, \quad L_{(1)}^M:=\frac{\rm d}{{\rm d} \lambda}(L_\lambda \cdot~)\big|_{\lambda=0}\in {\rm End}~M~ \text{is diagonalizable}.
\end{equation}
As before, we say that $m\in M$ has {\it conformal weight} $\Delta(m)\in \mathbb{F}$ if it is an eigenvector of $L_{(1)}^M$ of eigenvalue $\Delta(m)$.}
\end{definition}

The presence of a Virasoro element is quite common in the context of LCAs. In Section \ref{sec4}, we will discuss several important examples with finite rank including the $\mathcal{W}(b)$, Schrödinger-Virasoro, and extended Schrödinger-Virasoro LCAs. Here we present  some infinite rank examples:
\begin{example} (\kern-2.5pt\cite{SXY2,BDK,WCY}) The following are conformal LCAs of infinite rank:
\begin{itemize}
  \item[(1)]  The loop Virasoro Lie conformal algebra $\mathcal{L}$ is a free $\mathbb{F}[\partial]$-module with basis $\{L_i| i\in\mathbb{Z}\}$, endowed with the $\lambda$-bracket $$[L_{i\lambda} L_j]=(-\partial-2\lambda)L_{i+j}, ~~\text{for}~~ i,j\in \mathbb{Z}.$$
It is a conformal LCA with Virasoro element $-L_0$.
  \item[(2)]  The general Lie conformal algebra $gc_N$ is the free $\mathbb{F}[\partial]$-module generated by $\{J_A^n| n\in\mathbb{Z}_+,A\in gl_N\}$, with $\lambda$-bracket 
  $$[J_{A\lambda}^m J_B^n]=\sum_{s=0}^{m} \dbinom{m}{s}(\lambda+\partial)^s J_{AB}^{m+n-s}- \sum_{s=0}^{n} \dbinom{n}{s}(-\lambda)^s J_{BA}^{m+n-s}, ~~\text{for}~~ m,n\in \mathbb{Z}_+.$$
The element $J_I^1$ is a Virasoro element of $gc_N$.
  \item[(3)]  The Lie conformal algebra of Block type $\mathfrak{B}(p)$ (with $p\in\mathbb{C}$) is a free $\mathbb{F}[\partial]$-module with basis $\{L_i| i\in\mathbb{Z}_+\}$ and $\lambda$-bracket $$[L_{i\lambda} L_j]=
  \big((i+p)\partial+(i+j+2p)\lambda\big)L_{i+j},~~\text{for}~~ i,j\in \mathbb{Z}_+. $$
  For $p\neq0$, the element $L_0/p$ is a Virasoro element.
\end{itemize}
\end{example}

Throughout the remainder of this section, we always assume that $\mathcal{A}$ is a free conformal LCA and $M$ is a free conformal $\mathcal{A}$-module with respect to the Virasoro element $L$ in $\mathcal{A}$, so that both $\mathcal{A}$ and $M$ have an $\mathbb{F}[\partial]$-basis consisting of eigenvectors of $L_{(1)}$ and $L_{(1)}^M$, respectively. By skew-symmetry, a $q$-cochain $f$ is determined by its value on $a_1\otimes\cdots \otimes a_q$, where $a_i$ is an element in the basis consisting of the $L_{(1)}$-eigenvectors. 


We extend the notion of conformal weight to the spaces of polynomials $M[\lambda]$ by setting $\Delta(m\lambda^n)=\Delta(m)+n$, i.e. we assign to $\lambda$ conformal weight $1$ and extend in the obvious way. In other words, the conformal weights in $M[\lambda]$ are the eigenvalues of the operator  
\begin{align}\label{energy operator}
 E &=L_{(1)}^M+\lambda \frac{\rm d}{{\rm d}\lambda}. 
\end{align}

We have the following lemma by a verbatim repetition of the proof of \cite[Lemma 3.23]{BDK}. 
\begin{lemma}
    Let $a\in \mathcal{A}$ and $m\in M$ have conformal weights $\Delta(a)$ and $\Delta(m)$, respectively. Then 
\begin{enumerate}
    \item[(1)] $E$ is diagonalizable on $M[\lambda]$;
    \item[(2)] $\Delta(\partial a)=\Delta(a)+1,~\Delta(\partial m)=\Delta(m)+1$;
    \item[(3)] $\Delta(a_\lambda m)=\Delta(a)+\Delta(m)-1$.
\end{enumerate}
\end{lemma}

More generally, the linear operator 
\begin{align}\label{energy operator2} E=L_{(1)}^M+\sum_{i=1}^{n}\lambda_i \frac{\rm d}{{\rm d}\lambda_i} \end{align}
is diagonalizable in $M[\lambda_1,\cdots,\lambda_n]$. We will call the operator $E$ given by \eqref{energy operator2} the {\it energy operator}, and its eigenvalues will be called {\it conformal weights.} If $f$ is a $q$-cochain in $\widetilde{C}^q(\mathcal{A},M)$, then $f_{\lambda_1,...,\lambda_q}(a_1,...,a_q)\in M[\lambda_1,...,\lambda_q]$, for $a_1,\cdots,a_q\in \mathcal{A}$. We say that $f_{\lambda_1,...,\lambda_q}(a_1,...,a_q)$ is an eigenvector of $E$ with conformal weight $\Delta\big(f_{\lambda_1,...,\lambda_q}(a_1,...,a_q)\big)$ if and only if  \[ E\big(f_{\lambda_1,...,\lambda_q}(a_1,...,a_q)\big)=\Delta\big(f_{\lambda_1,...,\lambda_q}(a_1,...,a_q)\big) f_{\lambda_1,...,\lambda_n}(a_1,...,a_q). \]  

Finally, we define a linear operator $\widetilde{E}$ on the basic complex $\widetilde{C}^q(\mathcal{A},M)$ by 
\begin{equation} \label{3.1}
(\widetilde{E}f)_{\lambda_1,...,\lambda_q}(a_1,...,a_q)=E\big(f_{\lambda_1,...,\lambda_q}(a_1,...,a_q)\big)+\big(q-\sum_{i=1}^{q} \Delta(a_i)\big)f_{\lambda_1,...,\lambda_q}(a_1,...,a_q),
\end{equation} 
which we will call again the {\it energy operator.} As before, we call {\it conformal weights} the eigenvalues of the energy operator $\widetilde{E}$ in \eqref{3.1},
and we denote by $\Delta(f)$ the eigenvalue of the eigenvector $f\in\widetilde{C}^q(\mathcal{A},M)$. 
By \eqref{3.1}, we have
\begin{equation} \label{3.2}
\Delta\big(f_{\lambda_1,...,\lambda_q}(a_1,...,a_q)\big)=\Delta(f)+\sum_{i=1}^{q}\Delta(a_i)-q.
\end{equation}
Therefore,  $f\in\widetilde{C}^q(\mathcal{A},M)$ is an eigenvector of $\widetilde{E}$ if and only if every nonzero $f_{\lambda_1,...,\lambda_q}(a_1,...,a_q)$ is an eigenvector of $E$.

For given $r\in\mathbb{F}$ and $f\in\widetilde{C}^q(\mathcal{A},M)$, we define an $\mathbb{F}$-linear map $$f^{(r)}:\mathcal{A}^{\otimes q}\longrightarrow M[\lambda_1,...,\lambda_q]$$ satisfying conformal antilinearity, such that
\[f^{(r)}(a_1\otimes...\otimes a_q)=f^{(r)}_{\lambda_1,...,\lambda_q}(a_1,...,a_q),\] 
where $f^{(r)}_{\lambda_1,...,\lambda_q}(a_1,...,a_q)$ is the sum of all $E$-eigenvectors in $f_{\lambda_1,...,\lambda_q}(a_1,...,a_q)$ with eigenvalue $$r+\sum_{i=1}^{q}\Delta(a_i)-q. $$
By \eqref{3.2}, it is straightforward to see that $f^{(r)}\in \widetilde{C}^q(\mathcal{A},M)$ is an $\widetilde{E}$-eigenvector with conformal weight $r$, and
\begin{align}\label{decompsiton}
 f=\sum_{r\in\mathbb{F}} f^{(r)}.
\end{align}
This implies that $\widetilde{E}$ is diagonalizable on $\widetilde{C}^q(\mathcal{A},M)$. 
\begin{remark}
  In {\rm \cite{Su}}, a $q$-cochain of the form $f^{(r)}$ is referred to as a homogenous $q$-cochain of degree $r$. We observe that $f^{(r)}$ is in fact an $\widetilde{E}$-eigenvector with conformal weight $r$. 
\end{remark}



\begin{lemma}
If $f\in\widetilde{C}^q(\mathcal{A},M)$ is an $\widetilde{E}$-eigenvector, then its differential ${\bf d} f\in\widetilde{C}^{q+1}(\mathcal{A},M)$ is also an $\widetilde{E}$-eigenvector and their conformal weights satisfy \[ \Delta({\bf d} f)=\Delta(f). \] 
Moreover, the operator $\widetilde{E}$ commutes with ${\bf d}$ in $\widetilde{C}^{\bullet}(\mathcal{A},M)$, i.e., $\widetilde{E}\circ {\bf d}={\bf d}\circ \widetilde{E}$.  
\end{lemma}
\begin{proof} By \eqref{2.2},
we compute conformal weight for each term in the expansion of $({\bf d}f)_{\lambda_1,...,\lambda_{q+1}}(a_1,...,a_{q+1})$:
\begin{align*}
&\Delta\big(a_{i\lambda_i}f_{\lambda_1,...,\hat{\lambda}_i,...,\lambda_{q+1}}(a_1,...,\hat{a}_i,...,a_{q+1})\big)\\ 
&\quad\quad\quad =\Delta(a_i)+\Delta\big(f_{\lambda_1,...,\hat{\lambda}_i,...,\lambda_{q+1}}(a_1,...,\hat{a}_i,...,a_{q+1})\big)-1\\
&\quad\quad\quad =\Delta(a_i)+\Delta(f)+\sum_{j=1;j\neq i}^{q+1}\Delta(a_j)-q-1\\
&\quad \quad\quad =\Delta(f)+\sum_{i=1}^{q+1}\Delta(a_i)-(q+1),\\
&\Delta\big(f_{\lambda_i+\lambda_j,\lambda_1,...,\hat{\lambda}_i,...,\hat{\lambda}_j,...,\lambda_{q+1}}([a_{i\lambda_i}a_j],a_1,...,\hat{a}_i,...,\hat{a}_j,...,a_{q+1})\big)\\ &\quad\quad\quad=\Delta(f)+\Delta([a_{i\lambda_i}a_j])+\sum_{k=1;k\neq i,j}^{q+1}\Delta(a_k)-q \\
&\quad\quad\quad=\Delta(f)+\sum_{i=1}^{q+1}\Delta(a_i)-(q+1).
\end{align*}
Hence, $({\bf d}f)_{\lambda_1,...,\lambda_{q+1}}(a_1,...,a_{q+1})$ is an $E$-eigenvector of conformal weights $\Delta(f)+\sum_{i=1}^{q+1}\Delta(a_i)-(q+1).$ By \eqref{3.2}, we have $\Delta({\bf d}f)=\Delta(f). $
 Finally, since every cochain $f$ has a decomposition $f=\sum_{r\in\mathbb{F}} f^{(r)}$ , we have 
\[ \widetilde{E}\circ{\bf d}f=\sum_{r\in\mathbb{F}} \widetilde{E}\circ({\bf d}f^{(r)})=\sum_{r\in\mathbb{F}} r({\bf d}f^{(r)})={\bf d}\sum_{r\in\mathbb{F}} r f^{(r)}={\bf d}\circ \widetilde{E} f, \]
and thus $\widetilde{E}\circ {\bf d}={\bf d}\circ \widetilde{E}.$
\end{proof}
The following result will play an important role in the computation of basic cohomology.
\begin{lemma} \label{3.2.3} Let $f\in\widetilde{C}^q(\mathcal{A},M)$ with a decomposition $f=\sum_{r\in\mathbb{F}} f^{(r)}$. Then  
\[f=\sum_{r\in \mathbb{F}}f^{(r)}\in \widetilde{D}^q(\mathcal{A},M) \Longleftrightarrow f^{(r)}\in\widetilde{D}^q(\mathcal{A},M),~\forall r\in\mathbb{F}; \]
\[f=\sum_{r\in \mathbb{F}}f^{(r)}\in \widetilde{B}^q(\mathcal{A},M) \Longleftrightarrow f^{(r)}\in\widetilde{B}^q(\mathcal{A},M),~\forall r\in\mathbb{F}. \]
\end{lemma}
\begin{proof}As ${\mathbf d} f=\sum_{r\in \mathbb{F}} {\mathbf d}f^{(r)}$ decomposes ${\mathbf d} f$ into the sum of $\widetilde{E}$-eigenvectors of different conformal weights, and the terms cannot cancel if different from zero, we obtain  \[{\mathbf d} f=\sum_{r\in \mathbb{F}} {\mathbf d}f^{(r)}=0 \Longleftrightarrow {\mathbf d} f^{(r)}=0,~\forall r\in\mathbb{F}.\] 
Moreover, if there exists $\phi=\sum_{r\in \mathbb{F}} \phi^{(r)}\in \widetilde{C}^{q-1}$, such that $f={\mathbf d}\phi$, then 
\[ f={\mathbf d}\phi=\sum_{r\in \mathbb{F}}{\mathbf d}\phi^{(r)} \Longleftrightarrow f^{(r)}={\mathbf d}\phi^{(r)},~\forall r\in\mathbb{F}, \]
thus proving the lemma.
\end{proof}

Let us now consider an alternative interpretation of the energy operator $\widetilde{E}$. 

Define an operator $\tau:~\widetilde{C}^q(\mathcal{A},M)\longrightarrow\widetilde{C}^{q-1}(\mathcal{A},M)$ as follows: for $q=0$, set $\tau f=0$; for $q\geq1$, define
\begin{align*}
&(\tau f)_{\lambda_1,...,\lambda_{q-1}}(a_1,...,a_{q-1})=(-1)^{q-1}\frac{\rm d}{{\rm d}\lambda} f_{\lambda_1,...,\lambda_{q-1},\lambda}(a_1,...,a_{q-1},L)\big|_{\lambda=0}.
\end{align*}
For $q\geq 1$,  we have  
\begin{align*}
&\big(( {\bf d}\tau+\tau  {\bf d})f\big)_{\lambda_1,...,\lambda_q}(a_1,...,a_q) \\
=&~\sum_{k=1}^{q}(-1)^{k+1} a_{k\lambda_k}(\tau f)_{\lambda_1,...,\hat{\lambda}_k,...,\lambda_q}(a_1,...,\hat{a}_k,...,a_q) \\
&+\sum_{k,l=1;k<l}^{q}(-1)^{k+l} (\tau f)_{\lambda_k+\lambda_l,\lambda_1,...,\hat{\lambda}_k,...,\hat{\lambda}_l,...,\lambda_q}([a_{k\lambda_k}a_l],a_1,...,\hat{a}_k,...,\hat{a}_l,...,a_q) \\
&+(-1)^q \frac{\rm d}{{\rm d}\lambda} ({\mathbf d}f)_{\lambda_1,...,\lambda_q,\lambda}(a_1,...,a_q,L)\big|_{\lambda=0} \\
=&~(-1)^{q-1}\sum_{k=1}^{q}(-1)^{k+1}a_{k\lambda_k} \Big(\frac{\rm d}{{\rm d}\lambda}f_{\lambda_1,...,\hat{\lambda}_k,...,\lambda_q,\lambda}(a_1,...,\hat{a}_k,...,a_q,L)\big|_{\lambda=0}\Big) \\
&+(-1)^{q-1}\sum_{k,l=1;k<l}^{q}(-1)^{k+l} \frac{\rm d}{{\rm d}\lambda} f_{\lambda_k+\lambda_l,\lambda_1,...,\hat{\lambda}_k,...,\hat{\lambda}_l,...,\lambda_q,\lambda}([a_{k\lambda_k}a_l],a_1,...,\hat{a}_k,...,\hat{a}_l,...,a_q,L)\big|_{\lambda=0} \\
&+(-1)^q\sum_{k=1}^{q}(-1)^{k+1} \frac{\rm d}{{\rm d}\lambda}\big(a_{k\lambda_k} f_{\lambda_1,...,\hat{\lambda}_k,...,\lambda_q,\lambda}(a_1,...,\hat{a}_k,...,a_q,L)\big)\big|_{\lambda=0} \\
&+(-1)^q(-1)^{q+2} \frac{\rm d}{{\rm d}\lambda}\big(L_\lambda f_{\lambda_1,...,\lambda_q}(a_1,...,a_q)\big)\big|_{\lambda=0} \\
&+(-1)^q\sum_{k,l=1;k<l}^{q}(-1)^{k+l} \frac{\rm d}{{\rm d}\lambda}f_{\lambda_k+\lambda_l,\lambda_1,...,\hat{\lambda}_k,...,\hat{\lambda}_l,...,\lambda_q,\lambda}([a_{k\lambda_k}a_l],a_1,...,\hat{a}_k,...,\hat{a}_l,...,a_q,L)\big|_{\lambda=0} \\
&+(-1)^q\sum_{k=1}^{q}(-1)^{k+q+1} \frac{\rm d}{{\rm d}\lambda}f_{\lambda_k+\lambda,\lambda_1,...,\hat{\lambda}_k,...,\lambda_q}([a_{k\lambda_k}L],a_1,...,\hat{a}_k,...,a_q)\big|_{\lambda=0} \\
\triangleq&~ (a)+(b)+(c)+(d)+(e)+(f).
\end{align*}
It's easy to see that $(b)+(e)=0$ and $(a)+(c)=0$. 
If $f$ is an $\widetilde E$-eigenvector, then we can assume that 
\begin{align}\label{sum}
 f_{\lambda_1,...,\lambda_q}(a_1,...,a_q) = \sum_{j,k} f_{jk}(\lambda_1,...,\lambda_q)\partial^j m_k,
\end{align}
where $m_k$'s form an $\mathbb{F}[\partial]$-basis of $M$ consisting of $L_{(1)}^M$-eigenvectors and each $f_{jk}(\lambda_1,...,\lambda_q)$ is a homogenous polynomial 
in $\lambda_1,...,\lambda_q$. Note that for a fixed $q$-tuple of elements $a_1,...,a_q$, only a finite number of the summands in \eqref{sum} will produce values different from zero. 
According to \eqref{2.4} and conformal sesquilinearity, we have 
\begin{align*}
(d)&=\sum_{j,k} \frac{\rm d}{{\rm d}\lambda} \big(L_{\lambda}\partial^j m_k f_{jk}(\lambda_1,...,\lambda_q)\big)\big|_{\lambda=0} \\
&=\sum_{j,k} \frac{\rm d}{{\rm d}\lambda} \big((\lambda+\partial)^j(\partial+\Delta(m_k)\lambda)m_k+(\lambda+\partial)^j\mathrm{O}(\lambda^2)\big)\big|_{\lambda=0} f_{jk}(\lambda_1,...,\lambda_q) \\
&=\sum_{j,k} \big(j+\Delta(m_k)\big)\partial^j m_k f_{jk}(\lambda_1,...,\lambda_q) \\
&=\sum_{j,k} \Delta (\partial^j m_k)\partial^j m_k f_{jk}(\lambda_1,...,\lambda_q). 
\end{align*}
By conformal sesquilinearity again and \eqref{2.3}, we have
\[ [a_{i\lambda_i}L]=-\partial a_i+(\partial+\lambda_i)\Delta(a_i)a_i-\mathrm{O}\big((\lambda_i+\partial)^2\big), \] 
and thus
\begin{align*}
(f)&=\sum_{i=1}^{q} \frac{\rm d}{{\rm d}\lambda} f_{\lambda_1,...,\lambda_i+\lambda,...,\lambda_q}(a_1,...,[a_{i\lambda_i}L],...,a_q)\big|_{\lambda=0} \\
&=\sum_{i=1}^{q} \frac{\rm d}{{\rm d}\lambda} \big(\lambda_i+\big(1-\Delta(a_i)\big)\lambda\big) f_{\lambda_1,...,\lambda_i+\lambda,...,\lambda_q}(a_1,...,a_i,...,a_q)\big|_{\lambda=0}+\sum_{i=1}^{q} \frac{\rm d}{{\rm d}\lambda} \mathrm{O}(\lambda^2)\big|_{\lambda=0} \\
&=\sum_{i=1}^{q} \big(1-\Delta(a_i)\big)f_{\lambda_1,...,\lambda_i+\lambda,...,\lambda_q}(a_1,...,a_i,...,a_q)\big|_{\lambda=0} +\sum_{i=1}^{q} \lambda_i \frac{\rm d}{{\rm d}\lambda} f_{\lambda_1,...,\lambda_i+\lambda,...,\lambda_q}(a_1,...,a_i,...,a_q)\big|_{\lambda=0} \\
&=\sum_{i=1}^{q}\big(1-\Delta(a_i)\big) f_{\lambda_1,...,\lambda_q}(a_1,...,a_q) +\sum_{j,k} \partial^j m_k \sum_{i=1}^{q}\frac{\rm d}{{\rm d}\lambda} \lambda_i f_{jk}(\lambda_1,...,\lambda_i+\lambda,...,\lambda_q)\big|_{\lambda=0} \\
&=\sum_{i=1}^{q}\big(1-\Delta(a_i)\big) f_{\lambda_1,...,\lambda_q}(a_1,...,a_q)+\sum_{j,k} \deg (f_{jk})\partial^j m_k  f_{jk}(\lambda_1,...,\lambda_q). 
\end{align*}
By \eqref{3.2}, we have 
\begin{align*}
\Delta(f)&=\Delta\big(f_{\lambda_1,...,\lambda_q}(a_1,...,a_q)\big)+q-\sum_{i=1}^{q}\Delta(a_i) =\Delta\big(f_{jk}(\lambda_1,...,\lambda_q)\partial^j m_k\big)+\sum_{i=1}^{q}\big(1-\Delta(a_i)\big) \\
         &=\deg (f_{jk})+\Delta(\partial^j m_k)+\sum_{i=1}^{q}\big(1-\Delta(a_i)\big).
\end{align*}
Therefore, we obtain 
\begin{align*} 
&\big(({\mathbf d}\tau+\tau {\mathbf d})f\big)_{\lambda_1,...,\lambda_q}(a_1,...,a_q)
=(d)+(f) \\
&=\sum_{i=1}^{q}\big(1-\Delta(a_i)\big) f_{\lambda_1,...,\lambda_q}(a_1,...,a_q)+\sum_{j,k} \big(\Delta (\partial^j m_k)+\deg (f_{jk})\big)\partial^j m_k  f_{jk}(\lambda_1,...,\lambda_q) \\
&=\sum_{j,k} \Big(\sum_{i=1}^{q}\big(1-\Delta(a_i)\big)+\Delta (\partial^j m_k)+\deg (f_{jk})\Big) \partial^j m_k  f_{jk}(\lambda_1,...,\lambda_q) \\
&=\Delta(f)f_{\lambda_1,...,\lambda_q}(a_1,...,a_q),
\end{align*}
that is, \begin{equation} \label{3.3}
({\mathbf d}\tau+\tau {\mathbf d})f=\Delta(f)f=\widetilde{E} f.
\end{equation}
This equation holds for all cochains $f \in \widetilde{C}^q(\mathcal{A},M)$ with $q \geq 1$, as $f$ admits a decomposition into $\widetilde{E}$-eigenvectors. For the case $q=0$ where $f \in M$, the equality \eqref{3.3} remains valid, as verified by the following direct computation:

\[ ({\mathbf d}\tau+\tau {\mathbf d})f=\tau {\mathbf d}f=\frac{\rm d}{{\rm d}\lambda} ({\bf d}f)_{\lambda}(L)\big|_{\lambda=0}=\frac{\rm d}{{\rm d}\lambda} L_\lambda f\big|_{\lambda=0}=L_{(1)}^M f=\widetilde{E} f. \]

Based on the above discussion, we now state the main theorem, which plays a fundamental role in computing the basic cohomology.
\begin{theorem} \label{3.2.4} Let $\mathcal{A}$ be a free LCA containing a Virasoro element $L$, and $M$ a free conformal $\mathcal{A}$-module with respect to $L$. If an $\widetilde{E}$-eigenvector $f\in\widetilde{C}^q(\mathcal{A},M)$ is a cocycle, but not a coboundary, then its conformal weight equals zero.
\end{theorem}
\begin{proof}
From \eqref{3.3} and ${\bf d}f=0$, we have ${\bf d}\tau f=\Delta(f)f$. If $\Delta(f)\neq0$, then $f\in \widetilde{B}^q(\mathcal{A},M)$ when $q\geq1$ and $f=0$ when $q=0$. This is a contradiction.
\end{proof}

In conclusion, to compute the basic cohomology, Lemma \ref{3.2.3} shows that we only need to consider the 
$\widetilde{E}$-eigenvectors, and Theorem \ref{3.2.4} further indicates that these eigenvectors have conformal weight zero.



Finally, we discuss the computation of reduced cohomology. Denote by $[f]$ the cohomology class of $f$ and let $\bar{f}=f+\langle\partial\rangle$, where $f\in \widetilde{C}^\bullet(\mathcal{A},M)$. 
According to \cite{BKV}, we have \[ \widetilde{\mathrm{H}}^q(\mathcal{A},M)\cong \mathrm{H}^q\big(\partial\widetilde{C}^\bullet(\mathcal{A},M)\big),~q\geq 1, \]
and the isomorphism is given by $[f]\longmapsto [\partial f]$. Moreover, the short exact sequence $0\rightarrow\partial\widetilde{C}^\bullet\xrightarrow{\iota}\widetilde{C}^\bullet\xrightarrow{\pi}C^\bullet\rightarrow 0$ gives the long exact sequence of the cohomology groups:
\begin{equation} \label{3.4}
\begin{split}
    \cdots &\longrightarrow ~\mathrm{H}^q\big(\partial\widetilde{C}^\bullet(\mathcal{A},M)\big)~~ \xrightarrow{~\iota_q~} ~\widetilde{\mathrm{H}}^q(\mathcal{A},M)~~ \xrightarrow{~\pi_q~} ~\mathrm{H}^q(\mathcal{A},M)~~ \xrightarrow{~\omega_q~} \\
    &\longrightarrow \mathrm{H}^{q+1}\big(\partial\widetilde{C}^\bullet(\mathcal{A},M)\big) \xrightarrow{\iota_{q+1}} \widetilde{\mathrm{H}}^{q+1}(\mathcal{A},M) \xrightarrow{\pi_{q+1}} \mathrm{H}^{q+1}(\mathcal{A},M) \xrightarrow{\omega_{q+1}} \cdots,
\end{split}
\end{equation}
where $\iota_q$ is the $q$-th natural embedding, $\pi_q$ is the $q$-th natural projection and $\omega_q$ is the $q$-th connecting homomorphism, which is given by $\omega_q([\bar{f}])=[{\bf d}f]$.  

For $[\partial f]\in \mathrm{H}^q\big(\partial\widetilde{C}^\bullet(\mathcal{A},M)\big)$, we have $\iota_q([\partial f])=[\partial f]\in \widetilde{\mathrm{H}}^q(\mathcal{A},M)$. However, the conformal weight of $\partial f$ in $\widetilde{C}^q(\mathcal{A},M)$ is $\Delta(f)+1$, which means that if $[f]\in\widetilde{\mathrm{H}}^q(\mathcal{A},M)\cong \mathrm{H}^q\big(\partial\widetilde{C}^\bullet(\mathcal{A},M)\big)$ is nontrivial, then $[\partial f]$ is trivial. Thus $\iota_q([\partial f])=[\partial f]=0$, and from the exact sequence in \eqref{3.4}, 
\[ \ker (\pi_q)=\mathrm{im}~ (\iota_q)=\{0\},\quad \mathrm{im}~(\omega_q)=\ker (\iota_{q+1})=\mathrm{H}^{q+1}\big(\partial\widetilde{C}^\bullet(\mathcal{A},M)\big). \]
Finally we obtain a short exact sequence:\[ 0\longrightarrow \widetilde{\mathrm{H}}^q(\mathcal{A},M) \xrightarrow{~\pi_q~} \mathrm{H}^q(\mathcal{A},M) \xrightarrow{~\omega_q~} \mathrm{H}^{q+1}\big(\partial\widetilde{C}^\bullet(\mathcal{A},M)\big) \longrightarrow 0,\] 
which gives that
\begin{align}
\dim~ \mathrm{H}^q(\mathcal{A},M)&=\dim~ \widetilde{\mathrm{H}}^q(\mathcal{A},M) + \dim~ \mathrm{H}^{q+1}\big(\partial\widetilde{C}^\bullet(\mathcal{A},M)\big) \nonumber\\
&=\dim~ \widetilde{\mathrm{H}}^q(\mathcal{A},M) + \dim~ \widetilde{\mathrm{H}}^{q+1}(\mathcal{A},M),\quad \forall~q\geq1.\label{dim formular}
\end{align}
And a basis of $\mathrm{H}^q(\mathcal{A},M)$ is the combination of the projections of a basis of $\widetilde{\mathrm{H}}^q(\mathcal{A},M)$ and the pre-images of a basis of $\widetilde{\mathrm{H}}^{q+1}(\mathcal{A},M)\cong \mathrm{H}^{q+1}\big(\partial\widetilde{C}^\bullet(\mathcal{A},M)\big)$. This is our way to get the reduced cohomology through the basic cohomology.

The pre-images can be computed as follows. Let $[g]\in\widetilde{\mathrm{H}}^{q+1}(\mathcal{A},M)$, i.e. $g$ is a nontrivial $(q+1)$-cocycle, which implies that $\Delta(g)=0$, and correspondingly $[\partial g]\in \mathrm{H}^{q+1}\big(\partial\widetilde{C}^\bullet(\mathcal{A},M)\big)\cong\widetilde{\mathrm{H}}^{q+1}(\mathcal{A},M)$, then
\[\omega_q\big([\overline{\tau(\partial g)}]\big)=[\mathbf d \big(\tau(\partial g)\big)]=[({\bf d}\tau+\tau {\bf d})(\partial g)]=[\widetilde{E}(\partial g)]=[(\Delta(g)+1)(\partial g)]=[\partial g]. \]
Therefore, the pre-image of $[\partial g]$ under the connecting homomorphism $\omega_q$ is $[\overline{\tau(\partial g)}].$

\section{The algorithm for computing cohomology of finite LCAs and its applications}\label{sec4}

In this section, we first present an algorithm to compute the basic and reduced cohomology for finitely generated Lie conformal algebras containing a Virasoro element, with coefficients in their conformal modules of finite rank. We then apply this algorithm to compute the cohomology of our main examples: the $\mathcal{W}(b)$ conformal algebra, the Schr\"odinger-Virasoro conformal algebra, and the extended Schr\"odinger-Virasoro conformal algebra, each with coefficients in the adjoint module.


\subsection{ Algorithm for computing cohomology of finite LCAs}

In this subsection, we always assume that $\mathcal{A}$ is a finitely free conformal LCA, and $M$ is a finitely free conformal $\mathcal{A}$-module with respect to the Virasoro element $L$ in $\mathcal{A}$. Furthermore, we assume that $\mathcal{A}$ and $M$ admit $\mathbb{F}[\partial]$-bases $\{a_1,...,a_n\}$ and $\{m_1,...,m_l\}$, respectively, where each basis consists of eigenvectors of $L_{(1)}$ and $L_{(1)}^M$.

Under these assumptions, each $q$-cochain corresponds bijectively to a finite-dimensional matrix whose entries are polynomials in $\mathbb{F}[\partial,\lambda_1,...,\lambda_q]$. The number of rows in the matrix depends on $q$ and the rank of $\mathcal{A}$, while the number of columns is equal to the rank of $M$. Furthermore, by skew-symmetry, we need only consider elements in $\mathcal{A}^{\otimes q}$ of the form:
\[\underbrace{a_1\otimes...\otimes a_1}_{s_1}\otimes \underbrace{a_2\otimes...\otimes a_2}_{s_2}\otimes...\otimes \underbrace{a_n\otimes...\otimes a_n}_{s_n},~~~s_1+...+s_n=q,~~s_i\geq0,~i=1,...,n.\] Thus, the number of rows in the matrix equals the number of nonnegative integer solutions to the equation $$s_1+...+s_n=q.$$ In this framework, a $q$-cochain is represented computationally as such a matrix. Furthermore, Lemma \ref{3.2.3} guarantees that each polynomial entry in the matrix is homogeneous, while Theorem \ref{3.2.4} precisely determines its degree. These structural constraints allow us to explicitly construct the matrix, including its coefficients.


As for the operator $\bf d$, we first construct the $\lambda$-actions of the basis elements $a_1,...,a_n$ on $M[\lambda_1,...,\lambda_q]$, which are transformations of vectors, and then we compute each different form of the terms in ${\bf d}f$ respectively. For example, $$a_{i\lambda_i}f_{\lambda_1,...,\hat{\lambda}_i,...,\lambda_{q+1}}(a_1,...,\hat{a}_i,...,a_{q+1})$$ is computed by inputting the corresponding row in $f$ to the function of $a_i$'s $\lambda_i$-action and $$f_{\lambda_i+\lambda_j,\lambda_1,...,\hat{\lambda}_i,...,\hat{\lambda}_j,...,\lambda_{q+1}}([a_{i\lambda_i}a_j],a_1,...,\hat{a}_i,...,\hat{a}_j,...,a_{q+1})$$ is computed by expanding the $\lambda_i$-bracket $[a_{i\lambda_i}a_j]$ and locating the corresponding row in $f$, which may require the skew-symmetry. And the function of $\bf d$ can be achieved.

After we have ${\bf d}f$, the rest part of the program is just generating and solving linear equations. The first system of linear equations to be solved is \[ \texttt{Coefficientlist}[{\bf d}f]=0. \]  Then we obtain a matrix $f$ that satisfies ${\bf d}f=0$ by substituting the solution. It's very likely that some coefficients still exist and we can construct a list of possible bases $\{f1,f2,...\}$, which are certain matrices without coefficients. Then we need to generate a $(q-1)$-cochain $g$, whose conformal weight is also zero due to the condition $f={\bf d}g$. Note that $g$ may have some coefficients. Solve in turn the linear equations constructed by \[ \texttt{Coefficientlist}[f1-{\bf d}g]=0,~\texttt{Coefficientlist}[f2-{\bf d}g]=0,... \]
We take the matrix that results in no solution, and then we get a new list more close to a basis. Finally, for each pair of matrices in the new list such as $f1$ and $f2$, we continue to construct and solve equations like \[ \texttt{Coefficientlist}[f1-Cf2-{\bf d}g]=0, \]
where $C$ is a scalar. This helps us delete some equivalent bases.


For the reduced cohomology, when $q\geq1$, we firstly compute the images of the basis of $\widetilde{\rm H}^{q+1}(\mathcal{A},M)$ under the operator $\tau\circ\partial$, which is just a transformation of matrices, and then combine them with the remainder of the basis of $\widetilde{\rm H}^q(\mathcal{A},M)$ modulo by $(\partial+\lambda_1+...,+\lambda_q)$. When $q=0$, since $C^0(\mathcal{A},M)=M/\partial M$, the basis consists of $M$'s bases whose images can be modulo by $(\partial+\lambda)$ under the $\lambda$-action of $\mathcal{A}$'s bases.

We can now determine both the basic and reduced $q$-th cohomology for any given $q$. While $q$ ranges over all positive integers ($q\in \mathbb{Z}_+$), making an infinite computation appear necessary, the following lemma demonstrates that only finitely many calculations are required.


\begin{lemma}
There exists an integer $N\geq0$, such that 
\[\widetilde{\rm H}^q(\mathcal{A},M)=0,~~{\rm H}^q(\mathcal{A},M)=0,\quad\forall~q>N.\] 
In other words, both the basic and reduced cohomology are of finite dimension.
\end{lemma}
\begin{proof} Assume that $\mathcal{A}$ has an $\mathbb{F}[\partial]$-basis $\{a_1,...,a_n\}$ and 
$M$ has an $\mathbb{F}[\partial]$-basis $\{m_1,...,m_l\}$. Furthermore, suppose that 
$a_i$ is an eigenvector of $L_{(1)}$ for $i=1,..., n$; and $m_j$ is an eigenvector of $L_{(1)}^M$ for $j=1,..., l$.

For each given $q\geq1$ and each basis element $\bar{a}$ in $\mathcal{A}^{\otimes q}$, denote by $s_i$ the number of $a_i$ in $\bar{a}$ for $i=1,...,n$, so that \[ s_1+...+s_n=q~~~~\text{and}~~~~s_i\geq0,~i=1,...,n. \]
According to the skew-symmetry, we can write $\bar{a}$ as \[\bar{a}=\underbrace{a_1\otimes...\otimes a_1}_{s_1}\otimes \underbrace{a_2\otimes...\otimes a_2}_{s_2}\otimes...\otimes \underbrace{a_n\otimes...\otimes a_n}_{s_n}.\] 
For an $\widetilde{E}$-eigenvector $f\in\widetilde{C}^q(\mathcal{A},M)$ of conformal weight zero, we have 
\[ \prod_{\underset{i<j}{i,j=1}}^{s_1}(\lambda_i-\lambda_j)\cdot\prod_{\underset{i<j}{i,j=s_1+1}}^{s_1+s_2}(\lambda_i-\lambda_j)\cdot~...~\cdot\prod_{\underset{i<j}{i,j=s_1+...+s_{n-1}+1}}^{q}(\lambda_i-\lambda_j)~\bigg|~f_{\lambda_1,...,\lambda_n}(\bar{a}). \]
Regarding $f_{\lambda_1,...,\lambda_n}(\bar{a})$ as a polynomial in $\lambda_1,...,\lambda_n$ with $\Delta(f)=0$, we have 
\begin{align*}
\sum_{i=1}^{n}\frac{s_i(s_i-1)}{2}\leq \deg \big(f_{\lambda_1,...,\lambda_n}(\bar{a})\big)&\leq \Delta\big(f_{\lambda_1,...,\lambda_n}(\bar{a})\big)-\min_{k=1,...,l} \Delta(m_k) 
=\sum_{i=1}^{n} s_i\Delta(a_i)-q-\min_{k=1,...,l} \Delta(m_k) \\
&\leq\max_{i=1,...,n} \Delta(a_i) \sum_{i=1}^{n} s_i-q-\min_{k=1,...,l} \Delta(m_k) 
=\max_{i=1,...,n} \Delta(a_i)q-q-\min_{k=1,...,l} \Delta(m_k).
\end{align*}
The second "$\leq$" arises from the homogeneous property and the decomposition of $f_{\lambda_1,...,\lambda_n}(\bar{a})$. Applying the mean inequality $$\frac{(s_1+...+s_n)^2}{n}\leq s_1^2+...+s_n^2,$$ we obtain
\[ \frac{q^2}{2n}-\frac{q}{2}\leq \sum_{i=1}^{n}\frac{s_i(s_i-1)}{2}\leq\max_{i=1,...,n} \Delta(a_i)q-q-\min_{k=1,..,l} \Delta(m_k).\]
If we denote $u=\max \{\Delta(a_i)|i=1,...,n\}$ and $v=\min\{\Delta(m_k)|k=1,...,l\}$, then the quadratic inequality becomes
\[ q^2+(1-2u)nq+2nv\leq0. \]
The discriminant is $\Delta_{n,u,v}=(n-2nu)^2-8nv$. If $\Delta_{n,u,v}<0$, then $\widetilde{\rm H}^q(\mathcal{A},M)=0$ for each $q\geq1$. Otherwise, we can also see that $\widetilde{\rm H}^q(\mathcal{A},M)=0$ when $q>(2nu-n+\sqrt{\Delta_{n,u,v}})/2$, i.e. $\exists~N_0\in\mathbb{Z}_+$, such that $\widetilde{\rm H}^q(\mathcal{A},M)=0$, $\forall ~q> N_0$.
Finally, using $\dim~ \mathrm{H}^q(\mathcal{A},M)=\dim~ \widetilde{\mathrm{H}}^q(\mathcal{A},M) + \dim~ \widetilde{\mathrm{H}}^{q+1}(\mathcal{A},M)$, we can get the conclusion.
\end{proof}
To summarize, we outline the algorithm as follows:
\begin{enumerate}
    \item[Step 1.] Input $q$ (and possibly other parameters).
    
    \item[Step 2.] Generate the matrix $f$ as described above.
    
    \item[Step 3.] Make $f$ skew-symmetric by constructing linear equations for the coefficients, solving them, and substituting the solutions back into $f$.
    
    \item[Step 4.] Construct the functions describing the $\lambda$-action of $a_i$ on $M[\lambda_1,\ldots,\lambda_q]$ and compute $\mathbf{d}f$.
    
    \item[Step 5.] Solve the linear system $\texttt{Coefficientlist}[\mathbf{d}f] = 0$ and substitute the solutions into $f$.
    
    \item[Step 6.] If $q \neq 0$ and $f \neq 0$, construct a list of potential bases $\{f_1, f_2, \ldots\}$. Then generate $g$ similarly to $f$ (but with different dimensions) and check for solutions to the linear systems
    \[\texttt{Coefficientlist}
    [f_1 - \mathbf{d}g] = 0,\quad \texttt{Coefficientlist}[f_2 - \mathbf{d}g] = 0,\ \ldots
    \]
    
    \item[Step 7.] Eliminate equivalent bases by solving
    \[
    \texttt{Coefficientlist}[f_1 - Cf_2 - \mathbf{d}g] = 0,
    \]
    where $f_1$, $f_2$ are distinct elements in the new list, to obtain a basis for $\widetilde{\mathrm{H}}^q(\mathcal{A},M)$.
    
    \item[Step 8.] For $q \geq 1$:
    \begin{itemize}
        \item Generate bases for $\widetilde{\mathrm{H}}^q(\mathcal{A},M)$ and $\widetilde{\mathrm{H}}^{q+1}(\mathcal{A},M)$;
        \item Divide the former by $(\partial + \lambda_1 + \cdots + \lambda_q)$;
        \item Apply the $\tau\circ\partial$ action (constructed via matrix transformations) to the latter.
    \end{itemize}
    For $q = 0$, verify whether each $a_{i\lambda}m_k$ is divisible by $(\partial + \lambda)$ for some $m_k$, thus computing the basis of $\mathrm{H}^q(\mathcal{A},M)$.        
\end{enumerate}


Note that the implementation requires prior knowledge of the specific $\lambda$-bracket of $\mathcal{A}$ and $\lambda$-action of $\mathcal{A}$ on $M$. While the programs differ for distinct choices of 
$\mathcal{A}$ and $M$, they share the same algorithmic framework and contain substantial common code components.

In general, for any finite and free Lie conformal algebra containing a Virasoro element, our algorithm applies to cases where the coefficients belong to: (1) the trivial module, (2) the adjoint module, and (3) nontrivial irreducible finitely free conformal modules. 


\begin{lemma}\label{lem1}
    Let $\mathcal{A}$ be a finite and free LCA containing a Virasoro element $L$, and let $M=\mathbb{F}[\partial]m$ be a free conformal $\mathcal{A}$-module of rank one.  Suppose that the $\lambda$-action of $L$ on $m$ is given by 
    \[ L_\lambda m=(\partial+\alpha+\Delta\lambda)m.\]
    If $\alpha\neq 0$, then $${\rm H}^\bullet(\mathcal{A},M)={\rm \widetilde{H}}^\bullet(\mathcal{A},M)=0.$$ 
\end{lemma}
\begin{proof}
     Let $f\in\widetilde{C}^q(\mathcal{A},M).$ Define an operator $\tau_0:\widetilde{C}^q(\mathcal{A},M)\longrightarrow\widetilde{C}^{q-1}(\mathcal{A},M)$ as follows: for $q=0,$ set $\tau_0f=0$; for $q\geq 1,$ define
    \[ (\tau_0f)_{\lambda_1,...,\lambda_{q-1}}(a_1,...,a_{q-1})=(-1)^{q-1}f_{\lambda_1,...,\lambda_{q-1},\lambda}(a_1,...,a_{q-1},L)\big|_{\lambda=0}. \]
    For $q\geq1$, a direct computation shows that
    \begin{align*}
    \big(({\bf d}\tau_0+\tau_0{\bf d})f\big)_{\lambda_1,...,\lambda_q}(a_1,...,a_q)=&L_{(0)}^Mf_{\lambda_1,...,\lambda_q}(a_1,...,a_q)-\sum_{i=1}^{q}f_{\lambda_1,...,\lambda_q}(a_1,...,L_{(0)}a_i,...,a_q) \\
    =&(\partial+\alpha)f_{\lambda_1,...,\lambda_q}(a_1,...,a_q)+\sum_{i=1}^q \lambda_i f_{\lambda_1,...,\lambda_q}(a_1,...,a_q) \\
    \equiv&\alpha f_{\lambda_1,...,\lambda_q}(a_1,...,a_q)~\big({\rm mod}~~\partial\widetilde{C}^q(\mathcal{A},M)\big).
    \end{align*}
   For $q=0,$ we have $({\bf d}\tau_0+\tau_0{\bf d})m=\tau_0{\bf d}m=({\bf d}m)_\lambda L\big|_{\lambda=0}=L_\lambda m\big|_{\lambda=0}=(\partial+\alpha)m\equiv \alpha m~({\rm mod}~~\partial M).$
    
    Observe that $\tau_0$ commutes with $\partial$. Thus, if ${\bf d}f\in\partial\widetilde{C}^{q+1}(\mathcal{A},M),$ then $\tau_0{\bf d}f\in\partial\widetilde{C}^{q}(\mathcal{A},M).$ Consequence, unless $\alpha=0,$ any cocycle $f$ is a reduced coboundary, implying ${\rm H}^\bullet(\mathcal{A},M)=0$. By \eqref{dim formular} and a direct computation on any zero-cocycle, it follows that ${\rm \widetilde{H}}^\bullet(\mathcal{A},M)=0.$  
\end{proof}


\subsection{Cohomology of the \texorpdfstring{$\mathcal{W}(b)$}{\mathcal{W}(b)} Lie conformal algebra}

For $b\in\mathbb{F}$, the $\mathcal{W}(b)$ Lie conformal algebra is a free $\mathbb{F}[\partial]$-module generated by $L$ and $H$ satisfying:
   \begin{equation} \label{4.2}
    [L_\lambda L]=(\partial+2\lambda)L,\quad[L_\lambda H]=\big(\partial+(1-b)\lambda\big)H, \quad[H_\lambda L]=\big(-b\partial+(1-b)\lambda\big)H, \quad[H_\lambda H]=0.  
    \end{equation}
This algebra was originally introduced in \cite{XY} and further studied in \cite{LY,WY}. In the case of $b=0$, it reduces to the Heisenberg-Virasoro conformal algebra, which was first introduced in \cite{SY} and whose cohomology with trivial coefficients was investigated in \cite{YW1}. Setting $b=-1$ yields the $\mathcal{W}(2,2)$ Lie conformal algebra, whose structure theory was studied in \cite{YW2}.


It's easy to see that the Lie conformal algebra $\mathcal{W}(b)$ and its adjoint module satisfy our algorithm's requirements. Using our computational program, we can compute both the basic and reduced cohomology for arbitrary $b\in \mathbb{F}$. The program takes $b$ as an input parameter, allowing us to determine these dimensions for any specific choice of $b$. However, it cannot produce a general formula for the dimension of the (basic) cohomology group of $\mathcal{W}(b)$. For conciseness, we present below the computational results only for the cases $b=-1,0,1$.


\begin{theorem} \label{w0}
For Lie conformal algebra $\mathcal{W}(0)$, we have
\begin{equation*} 
\dim~ \widetilde{\mathrm{H}}^q\big(\mathcal{W}(0),\mathcal{W}(0)\big)=\begin{cases} 1,~~if~q=1,3,\\ 2,~~if~q=2,\\ 0,~~otherwise,  \end{cases} \text{and}~ \dim~ \mathrm{H}^q\big(\mathcal{W}(0),\mathcal{W}(0)\big)=\begin{cases} 3,~~if~q=1,2,\\ 1,~~if~q=0,3,\\ 0,~~otherwise. \end{cases}
\end{equation*}
\end{theorem}

The theorem is obtained from the following results computed by \textit{Mathematica} programming:

\begin{longtblr}[caption={The Basic Cohomology of $\mathcal{W}(0)$ with Coefficients in the Adjoint Module}]{colspec={p{1.15cm} m{7.5cm} m{4cm}},rowhead=1,rowsep = 4pt,colsep = 12pt}
    \hline
    \textbf{Order} & \textbf{Basis} & $\widetilde{\rm H}^q\big(\mathcal{W}(0),\mathcal{W}(0)\big)$ \\ 
    \hline
    $q=0$ & $0$ & $0$ \\ 
    $q=1$ & $f_{11}:L\mapsto H$ & $\mathbb{F}[f_{11}]$ \\ 
    $q=2$ & $f_{21}:L\otimes H\mapsto H$ & $\mathbb{F}[f_{21}]\oplus\mathbb{F}[f_{22}]$ \\ 
    & $f_{22}:L\otimes L\mapsto H(\lambda_1-\lambda_2)$ & \\ 
    $q=3$ & $f_{31}:L\otimes L\otimes H\mapsto H(\lambda_1-\lambda_2)$ & $\mathbb{F}[f_{31}]$ \\ 
    $q\geq4$ & $0$ & $0$ \\
    \hline
\end{longtblr}

\begin{longtblr}[caption={The Reduced Cohomology of $\mathcal{W}(0)$ with Coefficients in the Adjoint Module},label={table2}]{colspec={p{1.15cm} m{7.5cm} m{4cm}},rowhead=1,rowsep = 4pt,colsep = 12pt}
    \hline
    \textbf{Order} & \textbf{Basis} & ${\rm H}^q\big(\mathcal{W}(0),\mathcal{W}(0)\big)$ \\ 
    \hline
    $q=0$ & $\bar{f}_{01}=H$ & $\mathbb{F}\int H$ \\ 
    $q=1$ & $\bar{f}_{11}:L\mapsto H$ & $\mathbb{F}[\bar{f}_{11}]\oplus\mathbb{F}[\bar{f}_{12}]\oplus\mathbb{F}[\bar{f}_{13}]$ \\ 
    & $\bar{f}_{12}:H\mapsto H$ &  \\
    & $\bar{f}_{13}:L\mapsto \partial H$ \\
    $q=2$ & $\bar{f}_{21}:L\otimes H\mapsto H$ & $\mathbb{F}[\bar{f}_{21}]\oplus\mathbb{F}[\bar{f}_{22}]\oplus\mathbb{F}[\bar{f}_{23}]$ \\ 
    & $\bar{f}_{22}:L\otimes L\mapsto H(\lambda_1-\lambda_2)$ & \\ 
    & $\bar{f}_{23}:L\otimes L\mapsto H(\partial+\lambda_2)$ & \\
    $q=3$ & $\bar{f}_{31}:L\otimes L\otimes H\mapsto H(\lambda_1-\lambda_2)$ & $\mathbb{F}[\bar{f}_{31}]$ \\ 
    $q\geq4$ & $0$ & $0$ \\
    \hline
\end{longtblr}

\begin{corollary}
\begin{enumerate}
\item[(1)] Every Casimir element of $\mathcal{W}(0)$ is a scalar multiple of $\int H.$
\item[(2)] Every derivation of $\mathcal{W}(0)$ is a linear combination of an inner derivation and derivations of the form  $\bar{f}_{11},$ $\bar{f}_{12}$ and $\bar{f}_{13}$ given in Table \ref{table2}.
\item[(3)] Up to equivalence, any first-order deformation of $\mathcal{W}(0)$, that preserves the product and the $\mathbb{F}[\partial]$-module structure, corresponds to a linear combination of the $2$-cocycles $\bar{f}_{21},~\bar{f}_{22},$ and $\bar{f}_{23}$ given in Table \ref{table2}.
\end{enumerate}
\end{corollary}


\begin{theorem} 
    For Lie conformal algebra $\mathcal{W}(1)$, we have
    \begin{equation*} 
    \dim~ \widetilde{\mathrm{H}}^q\big(\mathcal{W}(1),\mathcal{W}(1)\big)=\begin{cases} 2,~~if~q=2,3,\\ 0,~~otherwise,  \end{cases} \text{and}~  \dim~ \mathrm{H}^q\big(\mathcal{W}(1),\mathcal{W}(1)\big)=\begin{cases} 2,~~if~q=1,3,\\ 4,~~if~q=2,\\ 0,~~otherwise.  \end{cases}
    \end{equation*}
\end{theorem}
The theorem is obtained from the following results computed by \textit{Mathematica} programming:

\begin{longtblr}[caption={The Basic Cohomology of $\mathcal{W}(1)$ with Coefficients in the Adjoint Module}]{colspec={p{1.15cm} m{7.5cm} m{4cm}},rowhead=1,rowsep = 4pt,colsep = 12pt}
    \hline
    \textbf{Order} & \textbf{Basis} & $\widetilde{\rm H}^q\big(\mathcal{W}(1),\mathcal{W}(1)\big)$ \\ 
    \hline
    $q=0$ & $0$ & $0$ \\ 
    $q=1$ & $0$ & $0$ \\ 
    $q=2$ & $f_{21}:L\otimes H\mapsto H$ & $\mathbb{F}[f_{21}]\oplus\mathbb{F}[f_{22}]$ \\ 
    & $f_{22}:L\otimes L\mapsto H(\lambda_1^2-\lambda_2^2)$ & \\
    $q=3$ & $f_{31}:L\otimes L\otimes H\mapsto H(\lambda_1-\lambda_2)$ & $\mathbb{F}[f_{31}]\oplus\mathbb{F}[f_{32}]$ \\ 
    & $f_{32}:L\otimes L\otimes L\mapsto H(\lambda_1-\lambda_2)(\lambda_1-\lambda_3)(\lambda_2-\lambda_3)$  & \\
    $q\geq4$ & $0$ & $0$ \\
    \hline
\end{longtblr}

\begin{longtblr}[caption={The Reduced Cohomology of $\mathcal{W}(1)$ with Coefficients in the Adjoint Module},,label={table4}]{colspec={p{1.15cm} m{7.5cm} m{4cm}},rowhead=1,rowsep = 4pt,colsep = 12pt}
    \hline
    \textbf{Order} & \textbf{Basis} & ${\rm H}^q\big(\mathcal{W}(1),\mathcal{W}(1)\big)$ \\ 
    \hline
    $q=0$ & $0$ & $0$ \\ 
    $q=1$ & $\bar{f}_{11}:H\mapsto H$ & $\mathbb{F}[\bar{f}_{11}]\oplus\mathbb{F}[\bar{f}_{12}]$ \\ 
    & $\bar{f}_{12}:L\mapsto -H\lambda_1^2$ & \\
    $q=2$ & $\bar{f}_{21}:L\otimes H\mapsto H$ & $\mathbb{F}[\bar{f}_{21}]\oplus\mathbb{F}[\bar{f}_{22}]\oplus\mathbb{F}[\bar{f}_{23}]\oplus$ \\ 
    & $\bar{f}_{22}:L\otimes L\mapsto H(\lambda_1^2-\lambda_2^2)$ & $\mathbb{F}[\bar{f}_{24}]$ \\
    & $\bar{f}_{23}:L\otimes H\mapsto H(\partial+\lambda_2)$ & \\
    & $\bar{f}_{24}:L\otimes L\mapsto H\lambda_1\lambda_2(\lambda_1-\lambda_2)$ & \\
    $q=3$ & $\bar{f}_{31}:L\otimes L\otimes H\mapsto H(\lambda_1-\lambda_2)$ & $\mathbb{F}[\bar{f}_{31}]\oplus\mathbb{F}[\bar{f}_{32}]$ \\ 
    & $\bar{f}_{32}:L\otimes L\otimes L\mapsto H(\lambda_1-\lambda_2)(\lambda_1-\lambda_3)(\lambda_2-\lambda_3)$  & \\
    $q\geq4$ & $0$ & $0$ \\
    \hline
\end{longtblr}

\begin{corollary}
    \begin{enumerate}
    \item[(1)] Every Casimir element of $\mathcal{W}(1)$ is trivial.
    \item[(2)] Every derivation of $\mathcal{W}(1)$ is a linear combination of an inner derivation and derivations of the form $\bar{f}_{11},$  $\bar{f}_{12}$  given in Table \ref{table4}.
    \item[(3)] Up to equivalence, any first-order deformation of $\mathcal{W}(1)$, which preserves the product
and the $\mathbb{F}[\partial]$-module structure, corresponds to a linear combination of the $2$-cocycles $\bar{f}_{21},\bar{f}_{22},\bar{f}_{23}$ and $\bar{f}_{24}$ given in Table \ref{table4}.
    \end{enumerate}
\end{corollary}


\begin{theorem} 
    For Lie conformal algebra $\mathcal{W}(-1)$, we have
    \begin{equation*} 
    \dim~ \widetilde{\mathrm{H}}^q\big(\mathcal{W}(-1),\mathcal{W}(-1)\big)=\begin{cases} 1,~~if~q=2,4,\\ 2,~~if~q=3,\\ 0,~~otherwise,  \end{cases} \text{and}~  \dim~ \mathrm{H}^q\big(\mathcal{W}(-1),\mathcal{W}(-1)\big)=\begin{cases} 1,~~if~q=1,4,\\ 3,~~if~q=2,3,\\ 0,~~otherwise.  \end{cases}
    \end{equation*}
\end{theorem}
The theorem is obtained from the following results computed by \textit{Mathematica} programming:

\begin{longtblr}[caption={The Basic Cohomology of $\mathcal{W}(-1)$ with Coefficients in the Adjoint Module}]{colspec={p{1.15cm} m{7.5cm} m{4cm}},rowhead=1,rowsep = 4pt,colsep = 12pt}
    \hline
    \textbf{Order} & \textbf{Basis} & $\widetilde{\rm H}^q\big(\mathcal{W}(-1),\mathcal{W}(-1)\big)$ \\ 
    \hline
    $q=0$ & $0$ & $0$ \\ 
    $q=1$ & $0$ & $0$ \\ 
    $q=2$ & $f_{21}:L\otimes H\mapsto H$ & $\mathbb{F}[f_{21}]$ \\ 
    $q=3$ & $f_{31}:L\otimes H\otimes H\mapsto L(\lambda_2-\lambda_3)$ & $\mathbb{F}[f_{31}]\oplus\mathbb{F}[f_{32}]$ \\ 
    & $f_{32}:L\otimes L\otimes H\mapsto H(\lambda_1-\lambda_2)$  & \\
    $q=4$ & $f_{41}:L\otimes L\otimes H\otimes H\mapsto L(\lambda_1-\lambda_2)(\lambda_3-\lambda_4)$ & $\mathbb{F}[f_{41}]$ \\
    $q\geq5$ & $0$ & $0$ \\
    \hline
\end{longtblr}

\begin{longtblr}[caption={The Reduced Cohomology of $\mathcal{W}(-1)$ with Coefficients in the Adjoint Module},label={table6}]{colspec={p{1.15cm} m{7.5cm} m{4cm}},rowhead=1,rowsep = 4pt,colsep = 12pt}
    \hline
    \textbf{Order} & \textbf{Basis} & ${\rm H}^q\big(\mathcal{W}(-1),\mathcal{W}(-1)\big)$ \\ 
    \hline
    $q=0$ & $0$ & $0$ \\ 
    $q=1$ & $\bar{f}_{11}:H\mapsto H$ & $\mathbb{F}[\bar{f}_{11}]$ \\ 
    $q=2$ & $\bar{f}_{21}:L\otimes H\mapsto H$ & $\mathbb{F}[\bar{f}_{21}]\oplus\mathbb{F}[\bar{f}_{22}]\oplus\mathbb{F}[\bar{f}_{23}]$ \\ 
    & $\bar{f}_{22}:H\otimes H\mapsto L(\lambda_1-\lambda_2)$ & \\
    & $\bar{f}_{23}:L\otimes H\mapsto H(\partial+\lambda_2)$ & \\
    $q=3$ & $\bar{f}_{31}:L\otimes H\otimes H\mapsto L(\lambda_2-\lambda_3)$ & $\mathbb{F}[\bar{f}_{31}]\oplus\mathbb{F}[\bar{f}_{32}]\oplus\mathbb{F}[\bar{f}_{33}]$ \\ 
    & $\bar{f}_{32}:L\otimes L\otimes H\mapsto H(\lambda_1-\lambda_2)$  & \\
    & $\bar{f}_{33}:L\otimes H\otimes H\mapsto L(\lambda_2-\lambda_3)(-\lambda_1)$ & \\
    $q=4$ & $\bar{f}_{41}:L\otimes L\otimes H\otimes H\mapsto L(\lambda_1-\lambda_2)(\lambda_3-\lambda_4)$ & $\mathbb{F}[\bar{f}_{41}]$ \\
    $q\geq5$ & $0$ & $0$ \\
    \hline
\end{longtblr}

\begin{corollary}
    \begin{enumerate}
    \item[(1)] Every Casimir element of $\mathcal{W}(-1)$ is trivial.
    \item[(2)] Every derivation of $\mathcal{W}(-1)$ is a linear combination of an inner derivation and the derivation $\bar{f}_{11}$ given in Table \ref{table6}.
    \item[(3)] Up to equivalence, any first-order deformation of $\mathcal{W}(-1)$, which preserves the product
and the $\mathbb{F}[\partial]$-module structure, corresponds to a linear combination of the $2$-cocycles $\bar{f}_{21},\bar{f}_{22}$ and $\bar{f}_{23}$ given in Table \ref{table6}. 
    \end{enumerate}
\end{corollary}

We now consider the cohomology of $\mathcal{W}(b)$ with coefficients in the  nontrivial irreducible finite conformal modules.
\begin{theorem}{\rm (\kern-2.5pt\cite[Theorem 2.8]{WY})}
Any nontrivial irreducible finite $\mathcal{W}(b)$-module is a free $\mathbb{F}[\partial]$-module of rank one, defined as follows:   
\begin{enumerate}
\item[(1)] if $b\neq0$, then $ M_{\Delta,\alpha}=\mathbb{F}[\partial]m,~~L_\lambda m=(\partial+\alpha+\Delta\lambda)m,~~H_\lambda m=0,~~\text{where}~~\alpha,\Delta\in\mathbb{F},~\Delta\neq0; $
\item[(2)] if $b=0$, then $M_{\Delta,\beta,\alpha}=\mathbb{F}[\partial]m,~~L_\lambda m=(\partial+\alpha+\Delta\lambda)m,~~H_\lambda m=\beta m,~~\text{where}~~\alpha,\beta,\Delta\in\mathbb{F},~(\Delta,\beta)\neq(0,0). $
\end{enumerate}
\end{theorem}

By Lemma \ref{lem1}, it suffices to compute cohomology for $\mathcal{W}(b)$ with coefficients in $ M_{\Delta,0}$ or $M_{\Delta,\beta,0}$, 
 where our algorithm exactly applies. Actually, we can compute the basic and reduced cohomology for given values of $\Delta$ and $\beta$ by \textit{Mathematica}.  However, we cannot give a unified result for any $\Delta$ and $\beta$. So here we only focus on a specific case when $\beta\neq0$ and $b=0$.

Define an operator $\tau_1:\widetilde{C}^q(\mathcal{W}(0),M_{\Delta,\beta,0})\longrightarrow\widetilde{C}^{q-1}(\mathcal{W}(0),M_{\Delta,\beta,0})$ as follows: for $q=0$, set $\tau_1f=0$; for $q\geq 1,$ define 
\[ (\tau_1f)_{\lambda_1,...,\lambda_{q-1}}(a_1,...,a_{q-1})=(-1)^{q-1}f_{\lambda_1,...,\lambda_{q-1},\lambda}(a_1,...,a_{q-1},H)\big|_{\lambda=0}. \]
Following the same proof as in Lemma \ref{lem1}, we derive $$({\bf d}\tau_1+\tau_1{\bf d})f=\beta f.$$ This immediately yields the following vanishing result for cohomology:
\begin{theorem}\label{th1}
 If $\beta\neq0$, then $\widetilde{\rm H}^\bullet(\mathcal{W}(0),M_{\Delta,\beta,0})={\rm H}^\bullet(\mathcal{W}(0),M_{\Delta,\beta,0})=0.$
\end{theorem}

\subsection{Cohomology of the (extended) Schr\"odinger-Virasoro conformal algebra}

Recall the Schr\"odinger-Virasoro conformal algebra $\mathcal{SV}=\mathbb{F}[\partial]L\oplus\mathbb{F}[\partial]M\oplus\mathbb{F}[\partial]Y$ is a free $\mathbb{F}[\partial]$-module of rank 3, generated by $L$, $M$ and $Y$, and satisfying the following non-zero $\lambda$-brackets:
    \begin{align}
    &[L_\lambda L]=(\partial+2\lambda)L,\quad [L_\lambda M]=(\partial+\lambda)M,\quad [L_\lambda Y]=(\partial+\frac{3}{2}\lambda)Y, \label{sv1}\\
    &[M_\lambda L]=\lambda M,\quad\quad\quad [Y_\lambda L]=(\frac{1}{2}\partial+\frac{3}{2}\lambda)Y,~~~ [Y_\lambda Y]=(\partial+2\lambda)M. \label{sv2}
    \end{align}
The extended Schr\"odinger-Virasoro conformal algebra $\widetilde{\mathcal{SV}}=\mathbb{F}[\partial]L\oplus\mathbb{F}[\partial]M\oplus\mathbb{F}[\partial]Y\oplus\mathbb{F}[\partial]N$ is a free $\mathbb{F}[\partial]$-module of rank 4, equipped with a $\mathbb{F}[\partial]$-basis $\{L,M, Y,N\}$. Its non-zero 
$\lambda$-brackets are given by \eqref{sv1}--\eqref{sv2}, along with the following relations:
    \begin{align*}
    &[L_\lambda N]=(\partial+\lambda)N,\quad [N_\lambda L]=\lambda N,\quad [N_\lambda M]=2M, \\
    &[M_\lambda N]=-2M,\quad\quad~~ [N_\lambda Y]=Y,\quad\quad [Y_\lambda N]=-Y. 
    \end{align*}   
Obviously, the Schr\"odinger-Virasoro conformal algebra is a subalgebra of the extended Schr\"odinger-Virasoro conformal algebra. They were originally  introduced in \cite{SY}, constructed from the Schr\"odinger-Virasoro Lie algebra and the extended Schr\"odinger-Virasoro Lie algebra, respectively. Their module extensions were later studied in \cite{YL}. Recently, their cohomologies with trivial coefficients were computed in \cite{WL1} and \cite{WL2}, respectively. It is easy to see that both $\mathcal{SV}$ and $\widetilde{\mathcal{SV}}$ are conformal with  the Virasoro element $L$. To validate our algorithm, we recomputed trivial coefficient cohomologies for them, and our algorithm reproduced the same results as those in \cite{WL1} and \cite{WL2}. 

In what follows, we compute the basic and reduced cohomologies of the Schrödinger-Virasoro conformal algebras 
 $\mathcal{SV}$ and $\widetilde{\mathcal{SV}}$ 
  with coefficients in their adjoint modules, employing our algorithmic approach. The detailed computations, carried out using Mathematica programming, are summarized in Tables \ref{table7} and \ref{table8}.

\begin{theorem} \label{sv}
    For the Schr\"odinger-Virasoro conformal algebra $\mathcal{SV}$, we have
\begin{equation*} 
    \dim~ \widetilde{\mathrm{H}}^q(\mathcal{SV},\mathcal{SV})=\begin{cases} 1,~~if~q=1,3,\\ 2,~~if~q=2,\\ 0,~~otherwise,  \end{cases} \text{and}~  \dim~ \mathrm{H}^q(\mathcal{SV},\mathcal{SV})=\begin{cases} 1,~~if~q=0,3,\\ 3,~~if~q=1,2,\\ 0,~~otherwise.  \end{cases}
\end{equation*}
\end{theorem}


\begin{longtblr}[caption={The Basic Cohomology of $\mathcal{SV}$ with Coefficients in the Adjoint Module},label={table7}]{colspec={p{1.15cm} m{7.5cm} m{4cm}},rowhead=1,rowsep = 4pt,colsep = 12pt}
    \hline
    \textbf{Order} & \textbf{Basis} & $\widetilde{\rm H}^q(\mathcal{SV},\mathcal{SV})$ \\ 
    \hline
    $q=0$ & $0$ & $0$ \\ 
    $q=1$ & $f_{11}:L\mapsto M$ & $\mathbb{F}[f_{11}]$ \\ 
    $q=2$ & $f_{21}:L\otimes Y\mapsto Y/2,~L\otimes M\mapsto M$ & $\mathbb{F}[f_{21}]\oplus\mathbb{F}[f_{22}]$ \\ 
    & $f_{22}:L\otimes L\mapsto M(\lambda_1-\lambda_2)$ & \\ 
    $q=3$ & $f_{31}:L\otimes L\otimes Y\mapsto Y(\lambda_1-\lambda_2)/2$ & $\mathbb{F}[f_{31}]$ \\ 
    & $\quad\quad L\otimes L\otimes M\mapsto M(\lambda_1-\lambda_2)$ &  \\ 
    $q\geq4$ & $0$ & $0$ \\
    \hline
\end{longtblr}

\begin{longtblr}[caption={The Reduced Cohomology of $\mathcal{SV}$ with Coefficients in the Adjoint Module},label={table8}]{colspec={p{1.15cm} m{7.5cm} m{4cm}},rowhead=1,rowsep = 4pt,colsep = 12pt}
    \hline
    \textbf{Order} & \textbf{Basis} & ${\rm H}^q(\mathcal{SV},\mathcal{SV})$ \\ 
    \hline
    $q=0$ & $\bar{f}_{01}=M$ & $\mathbb{F}\int M$ \\ 
    $q=1$ & $\bar{f}_{11}:L\mapsto M$ & $\mathbb{F}[\bar{f}_{11}]\oplus\mathbb{F}[\bar{f}_{12}]\oplus\mathbb{F}[\bar{f}_{13}]$ \\
    & $\bar{f}_{12}:Y\mapsto Y/2,~M\mapsto M$ & \\
    & $\bar{f}_{13}:L\mapsto \partial M$ & \\
    $q=2$ & $\bar{f}_{21}:L\otimes Y\mapsto Y/2,~L\otimes M\mapsto M$ & $\mathbb{F}[\bar{f}_{21}]\oplus\mathbb{F}[\bar{f}_{22}]\oplus\mathbb{F}[\bar{f}_{23}]$ \\ 
    & $\bar{f}_{22}:L\otimes L\mapsto M(\lambda_1-\lambda_2)$ & \\ 
    & $\bar{f}_{23}:L\otimes Y\mapsto Y(\partial+\lambda_2)/2$ & \\
    & $\quad\quad L\otimes M\mapsto M(\partial+\lambda_2)$ & \\
    $q=3$ & $\bar{f}_{31}:L\otimes L\otimes Y\mapsto Y(\lambda_1-\lambda_2)/2$ & $\mathbb{F}[\bar{f}_{31}]$ \\ 
    & $\quad\quad L\otimes L\otimes M\mapsto M(\lambda_1-\lambda_2)$ & \\
    $q\geq4$ & $0$ & $0$ \\
    \hline
\end{longtblr}

\begin{corollary}
    \begin{enumerate}
    \item[(1)] Every Casimir element of $\mathcal{SV}$ is a scalar multiple of $\int M.$
    \item[(2)] Every derivation of $\mathcal{SV}$ is a linear combination of an inner derivation and the derivations $\bar{f}_{11},$ $\bar{f}_{12},$ $\bar{f}_{13}$ given in Table \ref{table8}.
    \item[(3)] Up to equivalence, any first-order deformation of $\mathcal{SV}$, which preserves the product
and the $\mathbb{F}[\partial]$-module structure, corresponds to a linear combination of the $2$-cocycles $\bar{f}_{21},\bar{f}_{22}$ and $\bar{f}_{23}$ given in Table \ref{table8}. 
    \end{enumerate}
\end{corollary}


The result is trivial for  $\widetilde{\mathcal{SV}}$, as stated in the following theorem:

\begin{theorem} 
    For the extended Schr\"odinger-Virasoro conformal algebra $\widetilde{\mathcal{SV}}$, we have
\begin{equation*} 
    \widetilde{\mathrm{H}}^\bullet(\widetilde{\mathcal{SV}},\widetilde{\mathcal{SV}})= \mathrm{H}^\bullet(\widetilde{\mathcal{SV}},\widetilde{\mathcal{SV}})=0.
\end{equation*}
As a result, every Casimir element of $\widetilde{\mathcal{SV}}$ is trivial and every derivation of $\widetilde{\mathcal{SV}}$ is inner.
\end{theorem}

Finally, we consider the case where the coefficients are in the nontrivial irreducible finite  conformal modules.

\begin{theorem}{\rm (\kern-2.5pt\cite{WY,YL})} Any nontrivial irreducible finite $\widetilde{\mathcal{SV}}$-module is a free $\mathbb{F}[\partial]$-module of rank one, denoted by $$  V_{\Delta,\beta,\alpha}=\mathbb{F}[\partial]m,$$ with the $\lambda$-action 
    defined by
\[L_\lambda m=(\partial+\alpha+\Delta\lambda)m,~~N_\lambda m=\beta m,~~M_\lambda m=Y_\lambda m=0,\]
where $\alpha,\beta,\Delta\in\mathbb{F},~(\Delta,\beta)\neq(0,0).$ 
\end{theorem}

By Lemma \ref{lem1}, we only need to consider the case $\alpha=0$. Moreover, we focus on the case $\beta\neq 0$. 
Define an operator $\tau_2:\widetilde{C}^q(\widetilde{\mathcal{SV}},V_{\Delta,\beta,0})\longrightarrow\widetilde{C}^{q-1}(\widetilde{\mathcal{SV}},V_{\Delta,\beta,0})$ as follows: if $q=0,~\tau_2f=0$; and if $q\geq 1,$
\[ (\tau_2f)_{\lambda_1,...,\lambda_{q-1}}(a_1,...,a_{q-1})=(-1)^{q-1}f_{\lambda_1,...,\lambda_{q-1},\lambda}(a_1,...,a_{q-1},N)\big|_{\lambda=0}. \]
 A straightforward calculation yields $$({\bf d}\tau_2+\tau_2{\bf d})f=(\beta-2s_2-s_3)f,$$ where $s_2$ and $s_3$ donote the number of $M$ and $Y$, respectively, in the tensor product 
 $a_1\otimes...\otimes a_q\in\widetilde{\mathcal{SV}}^{\otimes q}$. Consequently, we obtain the following result:

\begin{theorem}
If $\beta\notin \mathbb{Z}$ or $\beta\in\mathbb{Z}_{-}$, then $\widetilde{\rm H}^\bullet(\widetilde{\mathcal{SV}},V_{\Delta,\beta,0})={\rm H}^\bullet(\widetilde{\mathcal{SV}},V_{\Delta,\beta,0})=0$.
\end{theorem}

\end{document}